\documentclass[runningheads]{llncs}

\usepackage{graphicx}
\usepackage{amsmath,amssymb}
\usepackage{hyperref}
\usepackage[noabbrev,capitalize]{cleveref}
\usepackage{todonotes}
\usepackage{booktabs}
\usepackage{multirow}
\usepackage{algorithm}
\usepackage{algpseudocode}

\usepackage{todonotes}
\renewcommand{\leq}{\leqslant}
\renewcommand{\geq}{\geqslant}
\newcommand{\win}[1]{{\textbf{#1}}}

\begin{document}

\title{A Frank-Wolfe-based primal heuristic for quadratic mixed-integer optimization}
\titlerunning{A Frank-Wolfe-based heuristic for quadratic mixed-integer optimization}

\author{Gioni Mexi\inst{1,2} \and Deborah Hendrych\inst{1,2} \and Sébastien Designolle\inst{3} \and\\ Mathieu Besançon\inst{4} \and Sebastian Pokutta\inst{1,2}}
\authorrunning{Gioni Mexi et al.}

\institute{Zuse Institute Berlin, Germany \and Technische Universität Berlin, Germany \and Inria, ENS de Lyon, UCBL, LIP, France \and Université Grenoble Alpes, Inria, LIG, CNRS, France}

\maketitle

\begin{abstract}
We propose a primal heuristic for quadratic mixed-integer problems.
Our method extends the Boscia framework -- originally a mixed-integer
convex solver leveraging a Frank-Wolfe-based branch-and-bound approach -- to
address nonconvex quadratic objective functions and constraints.
We reformulate nonlinear constraints, introduce preprocessing steps, and a suite of
heuristics including rounding strategies, gradient-guided selection, and large
neighborhood search techniques that exploit integer-feasible vertices
generated during the Frank-Wolfe iterations.
Computational results demonstrate the effectiveness of our method in solving
challenging MIQCQPs, achieving improvements on QPLIB instances within minutes
and winning first place in the Land-Doig MIP Computational Competition 2025.
\end{abstract}
\section{Introduction}

Mixed-Integer Quadratically Constrained Quadratic Problems (MIQCQPs) represent a broad category in optimization,
capturing an array of applications in operations research and machine learning.
These problems combine the combinatorial difficulties of mixed-integer structures with nonlinear nonconvex constraints.
Despite this double difficulty, solution methods, algorithms, and solvers have been developed over the past decades for generic or specific forms of MIQCQPs, including Couenne~\cite{belotti2009couenne}, GloMIQO~\cite{misener2013glomiqo}, ANTIGONE~\cite{misener2014antigone}, BARON~\cite{zhang2024solving}, and SCIP~\cite{bestuzheva2023enabling,bolusani2024scipoptimizationsuite90,bestuzheva2023global}.
We refer the reader to~\cite{kronqvist2019review,kronqvist202550} for recent reviews of mixed-integer nonlinear optimization beyond quadratic functions.
Recently, commercial mixed-integer linear solvers such as Xpress and Gurobi added capabilities to solve MIQCQPs to global optimality, showing the strong interest from application areas.
We consider MIQCQPs in the following form:
\begin{align}
  \min_{x \in \mathbb{R}^n} \quad & \frac{1}{2} x^\top Q x + d^\top x \label{eq:miqcqp} \\
  \text{s.t.} \quad
  & x^\top A_k x + b_k^\top x + c_k \leq 0 &\quad\forall\,k=1,\dots,m, \nonumber\\
  & x^L_i \leq x_i \leq x^U_i &\quad\forall\,i\in N, \nonumber\\
  & x_i \in \mathbb{Z} &\quad\forall\,i\in I, \nonumber
\end{align}
where $n, m \in \mathbb{N}$ are the number of variables and constraints, $N := \{1,\dots,n\}$, $I \subseteq N$ is the index set of integer
variables, $Q, A_k \in \mathbb{Q}^{n \times n}$ are symmetric matrices, $d,\ b_k \in
\mathbb{Q}^n$, $c_k \in \mathbb{Q}$ for $k=1,\dots,m$, and $x^L, x^U \in
\bar{\mathbb{Q}}^n$ (with $\bar{\mathbb{Q}} := \mathbb{Q} \cup \{\pm\infty\}$) are the lower
and upper bounds of $x$. Note that both $Q$ and $A_k$ need not be positive semidefinite,
hence nonconvex quadratic problems are allowed. If $I = \emptyset$, the
problem reduces to a quadratically constrained quadratic program (QCQP).
In the following, when discussing a single quadratic constraint,
we may omit the constraint index $k$ for simplicity.

\begin{figure}[ht!]
  \centering
  \includegraphics[width=0.9\textwidth]{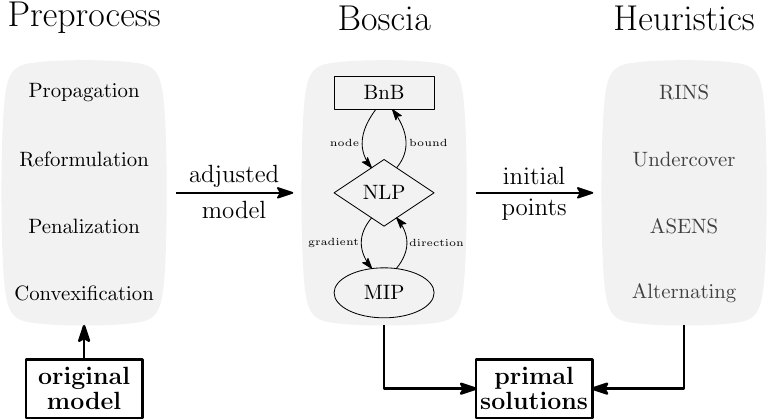}
  \caption{\label{fig:overview}
    \textbf{Overview of our approach.}
    The model first undergoes a pre-processing whose main goal is to transfer nonlinearities into the objective function.
    Our workhorse is indeed the Frank-Wolfe-based solver \emph{Boscia} which handles nonlinear problems (NLP) by suitably combining calls to a mixed-integer programming (MIP) solver in a branch-and-bound (BnB) framework.
    The points collected at each node in this process are then fed to various heuristics to reach better solutions.
  }
\end{figure}

In this work, we present a framework for this class of problems with an emphasis on finding primal feasible solutions, which we developed in the context of the Land-Doig MIP Computational Competition 2025~\footnote{Details on the competition can be found at the following page:\\ \url{https://www.mixedinteger.org/2025/competition/}}.
The framework also provides dual bounds for a subset of MIQCQPs as a by-product of the primal search.
Our approach is based on a branch-and-bound algorithm leveraging Frank-Wolfe to optimize continuous relaxations over the convex hull of mixed-integer feasible points,
as proposed in the Boscia framework~\cite{hendrych2022convex}.
Importantly, our method uses mixed-integer linear optimization solvers to compute vertices of a relaxed feasible region and can naturally handle differentiable nonlinear terms in the objective function, a property we leverage to transform and relax constraints.
Guided by the nonlinear relaxations, we fix continuous and integer variables and use other off-the-shelf solvers for the resulting smaller problems, adapting techniques from the large neighborhood search literature.
Our framework can thus be contrasted with most MINLP and MIQCQP solvers mentioned above which primarily construct polyhedral relaxations, with a few exceptions such as BARON.
A high-level overview of our framework is presented in Figure~\ref{fig:overview}.

The rest of the paper is organized as follows.
Section~\ref{sec:boscia} presents the core part of the framework consisting of the specialized branch-and-bound with relaxations solved with Frank-Wolfe.
Section~\ref{sec:problemtransformation} details the transformations we operate on the original MIQCQP, modifying the set of constraints and the objective function.
Section~\ref{sec:heuristics} expands on the heuristics called during and after the branch-and-bound process.
Section~\ref{sec:experiments} presents computational results assessing the performance of the framework and the impact of its different components.

\section{Boscia: Integer Frank-Wolfe}\label{sec:boscia}

Boscia~\cite{miconvexfw} is a mixed-integer convex solver that tackles problems of the form
\begin{align}
  \min_{x \in \mathbb{R}^n} \quad & f(x) \label{eq:boscia} \\
  \text{s.t.} \quad
  & x \in \mathcal{X} \nonumber\\
  & x_i \in \mathbb{Z} \quad\forall\,i\in I \nonumber
\end{align}
where $\mathcal{X}\subseteq \mathbb{R}^n$ is a convex, compact (polyhedral) set, $I\subseteq N$ is the set of integer variables and $f$ is a convex, continuously differentiable function and its gradient is Lipschitz continuous with constant $L$.

It operates a branch-and-bound scheme, utilizing Frank-Wolfe (FW) methods~\cite{braun2022conditional} implemented in the FrankWolfe.jl library~\cite{besanccon2022frankwolfe,besancon2025improvedalgorithmsnovelapplications} to solve continuous subproblems.
FW methods are projection-free first-order methods that only require oracle access to the objective function, its gradient, and the feasible region.
At iteration $t$, the classic FW algorithm, dubbed Vanilla Frank-Wolfe, calls the Linear Minimization Oracle (LMO) over the feasible region using the gradient of the current iterate $x^t$ as the objective function.
The solution $v^t$ is used to generate the next descent direction $d^t = x^t - v^t$, and the iterate is updated as $x^{t+1} = x^t + \gamma^t d^t$ where $\gamma^t$ is the step size.
Thus, the iterate can be represented as a convex combination of the extreme points of the feasible region.
This combination is referred to as the \emph{active set} in the FW context.
Many FW variants explicitly use the active set to make progress.
The Corrective Frank-Wolfe (CFW) algorithm~\cite{halbey2025efficient} shown in Algorithms~\ref{alg:cfw} and~\ref{alg:corrective-step} presents a framework for all active set based variants of FW.

In lines~\ref{alg:cfw:a_t}--\ref{alg:cfw:v_t}, the so-called \emph{away vertex} $a^t$, the vertex in the active set with the largest inner product with the gradient,
the \emph{local FW} $s^t$, vertex of the active set with the smallest inner product with the gradient, and the \emph{global FW} vertex $v^t$ are computed.
In line~\ref{alg:cfw:local_gap}, $a^t$ and $s^t$ are used to compute the local gap.
If it is sufficiently large, the algorithm proceeds to optimize the active set using the Corrective Step algorithm~\ref{alg:corrective-step}.
Note that the constant $L$ is the Lipschitz constant of the gradient of the objective function and $D$ is the diameter of the feasible region $\mathcal{X}$.
The specific corrective step depends on the variant of FW used.
The default variant used in Boscia, and in our set-up, is the Blended Pairwise Conditional Gradient (BPCG) method~\cite{TTP2021}.
The corrective step in this case shifts weight from the away vertex $a^t$ to the local FW vertex $s^t$, computed via linesearch.
If all of the weight is shifted to the local FW vertex, the away vertex is dropped (drop step), otherwise the descent step is performed with descent direction $d=a^t-s^t$.
As step-size rule, we use the secant line-search step-size rule~\cite{HBMP2025}, which is particularly efficient in the context of quadratic problems.
\begin{algorithm}[H]
  \caption{Corrective Frank-Wolfe (CFW) \cite{halbey2025efficient}}
  \label{alg:cfw}
  \begin{algorithmic}[1]
  \Require convex, smooth function $f$, start point $x^0 \in V(\mathcal{X})$ (vertex of $\mathcal{X}$).
  \State $S^0 \gets \{x^0\}$ \Comment{active set}
  \For{$t = 0$ to $T-1$}
    \State $a^t \gets \arg\max_{v \in S^t} \langle \nabla f(x^t), v \rangle$ \Comment{away vertex} \label{alg:cfw:a_t}
    \State $s^t \gets \arg\min_{v \in S^t} \langle \nabla f(x^t), v \rangle$ \Comment{local FW} \label{alg:cfw:s_t}
    \State $v^t \gets \arg\min_{v \in V(\mathcal{X})} \langle \nabla f(x^t), v \rangle$ \Comment{global FW} \label{alg:cfw:v_t}
    \If{$\langle \nabla f(x^t), a^t - s^t \rangle \geq \langle \nabla f(x^t), x^t - v^t \rangle$} \label{alg:cfw:local_gap}
      \State $x^{t+1}, S^{t+1} \gets \textsc{CorrectiveStep}(S^t, x^t, a^t, s^t)$
    \Else
      \State $\gamma^t \gets \arg\min_{\gamma \in [0,1]} f(x^t - \gamma(x^t - v^t))$
      \State $x^{t+1} \gets x^t - \gamma^t (x^t - v^t)$
      \State $S^{t+1} \gets S^t \cup \{v^t\}$
    \EndIf
  \EndFor
  \end{algorithmic}
  \end{algorithm}
\vspace{-0.7em}
\begin{algorithm}[H]
  \caption{Corrective Step($S, x, a, s$) \cite{halbey2025efficient}}
  \label{alg:corrective-step}
  \begin{algorithmic}[1]
  \Require $S \subset \mathcal{X}, \; x, a, s \in \mathcal{X}$
  \Ensure $S' \subseteq S, \; x' \in \operatorname{conv}(S')$ satisfying
    \Statex \(\displaystyle f(x') \le f(x) \quad\text{and}\quad S' \subsetneq S \)
    \hfill \(\triangleright\) \text{drop step}
    \Statex \(\displaystyle
    f(x) - f(x') \;\ge\;
    \frac{\langle \nabla f(\mathbf{x}),\, a-s\rangle^2}{2 L D^2} \)
    \hfill \(\triangleright\) \text{descent step}
  \end{algorithmic}
  \end{algorithm}
Boscia solves a Mixed-Integer Programming (MIP) problem as the LMO.
By doing so, Boscia optimizes convex relaxations over the convex hull of mixed-integer feasible points instead of the continuous relaxation of the constraints, 
resulting in a much smaller branch-and-bound search tree and directly leveraging the MIP solver machinery, e.g., cutting planes, conflict analysis, and heuristics.
This adjustment enables Boscia to inherently sample integer-feasible solutions, effectively embedding a heuristic search within the optimization process.
Importantly, optimizing over the convex hull of integer-feasible solutions is performed by Boscia without an explicit algebraic description of this feasible region.
In order to alleviate the cost of the multiple MIP solves, Boscia integrates lazified FW variants which reuse vertices of the feasible region computed by the MIP solver 
as much as possible, which can be performed without losing convergence guarantees; see, e.g.,~\cite{braun2019lazifying} for details on this lazification.
For further lazification, Boscia holds on to vertices dropped by FW as they might become useful for nodes later down the branch.
Hence, any integer-feasible point is computed at most once.
The active set is used during branching to facilitate warm-starting of the children nodes by splitting the active set of their parent.
Since all vertices are integer-feasible, any vertex $v$ of the active set must satisfy either $v_i \leq \lfloor \bar{x}_i \rfloor$ or $v_i \geq \lceil \bar{x}_i \rceil$ 
where $\bar{x}$ denotes the fractional node solution and $i$ the branching index.
This property is used to split the active set into two parts during branching.

MIQCQPs do not directly fit into the convex mixed-integer optimization framework that Boscia is designed for, therefore several key modifications to the algorithm are needed.
Since Boscia assumes linear constraints, the quadratic constraints must be reformulated to comply with this structure.
We achieve this by employing a power penalty relaxation, ensuring that the solver can process these constraints while maintaining feasibility (see Section~\ref{sec:penalty}).

Another key challenge is Boscia's assumption of convexity, which does not always hold for the given MIQCQPs.
To address this, we implement several modifications.
First, we disable node pruning based on the lower bound in the branch-and-bound process, as the lower bound information may no longer be valid.
Second, we instruct Boscia to ignore the standard FW lower bound, which is not applicable in our setting.
Finally, we adjust Boscia's solution storage mechanism.
By default, it performs pre-sampling and retains only improving solutions, but this approach proved inadequate for quadratic constraints.
Instead, we modify the solver to track all solutions encountered during the process, ensuring that valuable candidates are not prematurely discarded.
Additionally, we added a callback mechanism that allows us to evaluate the solutions with respect to the original objective and
lets us discard solutions which do not satisfy the quadratic constraints.
When Boscia finds a new solution, we check if it satisfies the quadratic constraints and discard it if it does not.
Furthermore, the solution is only added to the list of solutions if it improves upon the current best solution.

Because we are primarily interested in finding solutions, we do not let Boscia solve to completion.
Instead, we provide a node and time limit to the solver.
If the node limit is reached and there is still time remaining, we \emph{restart} Boscia with a new starting point. 
We compute the new start point via the LMO by using either the gradient of the current iterate or a random direction.
In this way, we aim to explore as much of the feasible region as possible.

To test multiple problem modifications, we employ a \emph{parallelization} strategy
as it proved to be computationally inefficient to change problem parameters between restarts.
For the convexification schemes, for example, we tried shifting different percentages of the spectrum of the Hessian into the positive orthant,
see Section~\ref{sec:convexification} for details.
For problems with quadratic constraints, we tried different power penalty parameters, see Section~\ref{sec:penalty} for details.
Note that the restarts are performed in each thread.

These modifications allow Boscia to effectively handle the challenges posed by quadratic mixed-integer problems, making it a competitive approach for the Land-Doig MIP Computational Competition 2025.
Even though the focus of the competition is on the primal side, we highlight that the developed solution framework is also suited to derive high-quality relaxation bounds,
since the constraint penalization and many transformations are exact reformulations or produce relaxations of the original problem.

\section{Problem transformations and presolving}\label{sec:problemtransformation}

In this section, we present the problem transformations performed before executing the main algorithms.
These presolving steps are crucial to obtain good performance on several problem classes.

\subsection{Propagation}
\label{sec:propagation}
Initially, we apply a propagation step to tighten the bounds of the variables.
For linear constraints, we apply a simple activity-based bound strengthening
technique and some specialized propagators for several classes of linear
constraints~\cite{fischetti2009feasibility}. For quadratic constraints, we
experimented with McCormick linearizations~\cite{mccormick1976computability},
introducing auxiliary variables $z_{ij}$ to replace bilinear terms $x_i x_j$,
and constructing McCormick envelopes. However, we found that
incorporating these linearized constraints in the problem formulation, or even just propagating them,
showed little to no performance improvement in practice.

\subsection{Power penalty}
\label{sec:penalty}

We take advantage of the capabilities of the Frank-Wolfe approach to solve arbitrary nonlinear objectives to relax quadratic constraints.
In particular, we leverage the power penalty formulation introduced in~\cite{sharma2024network} to relax quadratic constraints.
Consider an objective function $f(\cdot)$ and a quadratic constraint
\begin{align*}
  x^\top A x + b^\top x + c \leq 0.
\end{align*}
We reformulate the problem by dropping this constraint and integrating it into the objective which becomes
\begin{align*}
  f(x) + \mu \cdot {\max\{x^\top A x + b^\top x + c, 0\}}^p
\end{align*}
for a penalty parameter $\mu > 0$ and a real exponent $p > 1$.
A parameter $p=1$ would result in a nonsmooth, nondifferentiable function, while $p=2$ would result in a term akin to that of augmented Lagrangian methods.
We experimented with $p\in[1.2,1.8]$, trading off the steepness of the gradient towards feasible points and the smoothness of the continuous subproblems tackled by Frank-Wolfe.

\subsection{Quadratic special structures}

In addition to the power penalty reformulation, we also exploit the specific structure of two types of quadratic constraints: complementarities and perspective constraints.

\paragraph{Complementarities.} They are constraints of the form $x_i x_j = 0$.
Although they are quadratic equalities, they are combinatorial in nature and can be treated as such with an additional binary variable $z$ controlling which of the two terms should be set to zero.
In order to avoid relying on variable bounds for a big-M constraint, we reformulate the complementarity to a pair of indicator constraints, which we present here with the assumption that both variables are nonnegative:
\begin{align*}
  & z = 0 \Rightarrow x_i \leq 0, \\
  & z = 1 \Rightarrow x_j \leq 0.
\end{align*}

Importantly, we do not reformulate complementarity quadratic constraints to Special Ordered Sets of type 1 (SOS1) constraints, since such constraints would not respect the assumptions from Boscia (see~\cite{miconvexfw}).

\paragraph{Perspective constraints.} They are another type of special quadratic constraint from three variables: $x$ is a single variable, $w$ is the epigraph variable, and $z$ is a binary variable.
The constraint is of the form:
\begin{align*}
  & x^2 \leq z w\quad\text{where}\quad x \geq 0, w \geq 0, z \in \{0,1\}.
\end{align*}
A perspective constraint forms the convex hull of the set:
\begin{align*}
  \left\{(x,w,z) \in \mathbb{R}_+ \times \mathbb{R}_+ \times \{0,1\}\, :\, x^2 \leq w, (1-z) x =0 \right\},
\end{align*}
see, e.g., the seminal paper~\cite{frangioni2006perspective} on perspective functions in mixed-integer convex optimization and~\cite{furman2020computationally} and references therein on handling perspective functions in a nonlinear solver.
The binary variable activates the continuous variable $x$; the epigraph variable $w$ is often a consequence of a nonlinear objective term converted to a constraint.
When this is the case and the epigraph $w$ indeed only appears in the perspective constraint and in the objective, we transform the constraint back into a quadratic term $c x^2$ in the objective along with a big-M and/or indicator constraint explicitly tying $z$ and $x$.
This transformation separates the nonsmooth perspective constraint that can lead to numerical challenges in a nonlinear setting (as reported in~\cite{furman2020computationally} among others) into a smooth quadratic term in the objective and the combinatorial ``activation'' structure in the constraints.

\subsection{Convexification of quadratic binary problems}\label{sec:convexification}

A subclass of MIQCQPs are quadratic binary optimization problems, which involve objective functions of the form
\begin{equation}
  \min_{x \in \{0,1\}^n} \frac{1}{2}x^\top Q x + d^\top x,
\end{equation}
where $Q \in \mathbb{R}^{n \times n}$ is a symmetric matrix and $d \in \mathbb{R}^n$ is a linear term.
Since all variables are binary, we can use the identity $x_i^2 = x_i$ to reformulate the problem and ``increase'' convexity.

Formally, we replace $(Q, d)$ with $(Q_\ell,d_\ell)$, such that the convexified matrix $Q_\ell$ has a proportion $\ell$ of nonnegative eigenvalues, where $0 \leq \ell \leq 1$.
This is achieved through a spectral decomposition of $Q$, giving its eigenvalues $\lambda_1, \dots, \lambda_n$ from which we can easily obtain the index $i$ such that $Q_\ell=Q+\lambda_iI$ satisfies the desired property.
Accordingly, $d$ is mapped onto $d_\ell=d-\lambda_i/2$ to preserve the original objective value at all binary points.

This partial convexification balances problem tractability and the original problem structure, making convexification a tunable process controlled by the parameter $\ell$.
A whole stream of work studies hardness~\cite{pardalos1991quadratic} and algorithms~\cite{luo2019new} for quadratic optimization when a small number of eigenvalues are negative.
Given the finite time limit, there is indeed a trade-off between the deformation of the initial objective function, whose Lipschitz constant is strongly affected by the transformation, and the advantages of mitigating nonconvexity.
In practice, we explore different parameters $\ell$ in parallel on multiple threads.
Computational experiments (see Section~\ref{sec:convexification_results}) show that some convexification indeed helps but also confirm that full convexification can be detrimental to the search for primal solutions since the original objective is dominated by the convexifying term, sometimes by orders of magnitude.

Note that a much tighter convexification can be obtained via solving an SDP\@. To this end, we define a family of objective functions which are equivalent for all $x \in \{0,1\}^n$ via a parameter $u \in \mathbb{R}^n$:
\[
\frac{1}{2}x^\top (Q - 2 \operatorname{diag}(u)) x + (d+u)^\top x.
\]
The tightest possible convexification can then be obtained via the SDP that finds the optimal vector $u$ such that $Q - 2 \operatorname{diag}(u)$ is positive semidefinite while minimizing the objective value. This can be formulated as:
\begin{align*}
\min_{r, u}\;& r \\
\text{s.t.}\; &\begin{bmatrix} r & -(d + u)^\top / 2\\ -(d + u) / 2 & \operatorname{diag}(u) - Q/2 \end{bmatrix} \succeq 0.
\end{align*}
We did not explore this approach, though, as we were primarily interested in finding primal solutions and the tightness of the convexification was not critical to the empirical performance of our approach,
so that the additional computational cost of solving the SDP was not justified.

\section{Primal Heuristics}\label{sec:heuristics}
Our framework incorporates several primal heuristics to find feasible solutions from points generated during the Frank-Wolfe iterations.
\subsection{Rounding}
\label{sec:rounding_heuristics}
One of the simplest techniques is \emph{standard rounding}, where fractional values are rounded to the nearest integer.
Additionally, Boscia employs \emph{probability rounding} for binary variables, where variables are randomly fixed to 0 or 1 using the entries as probabilities, 
and Frank-Wolfe is then used to solve for the remaining continuous variables.
We present in Algorithm~\ref{alg:probability_rounding} the pseudo-code for this rounding technique.
Another specialized rounding method is designed for \emph{0/1 polytopes}, ensuring that the rounded solution remains within the feasible region.
Furthermore, for problems with \emph{simplex-like feasible regions}, Boscia implements a specialized rounding heuristic that 
is aware of the defining hyperplane structure, improving feasibility and efficiency in these cases.
Simplex-like feasible regions are the (scaled) unit simplex and the (scaled) probability simplex.
\begin{algorithm}[H]
\caption{Probability Rounding for a Mixed-Binary Program}
\label{alg:probability_rounding}
\begin{algorithmic}[1]
\Require fractional point $x$, set of integer variables $I$
\For{$i \in I$}
    \State $v \; \sim \; \text{Bernoulli}(x_i)$
    \If{$v = 1$}
        \State $\hat{x}_i \gets \min\{1, \lceil x_i \rceil\}$
    \Else
        \State $\hat{x}_i \gets \max\{0, \lfloor x_i \rfloor \} $
    \EndIf
\EndFor

\State $\hat{x} \gets$Fix bounds to the rounded entries and solve Frank-Wolfe for the continuous variables
\end{algorithmic}
\end{algorithm}

\subsection{Gradient-based heuristics}
Beyond the rounding strategies described above, we utilize \emph{gradient-based heuristics}, which explore a set of promising vertices, e.g., a follow-the-gradient heuristic inspired by the chasing gradients paradigm~\cite{CP2020boost}.
These heuristics are particularly suitable for Boscia since the algorithmic components they require (gradient and linear minimization oracles) are precisely those that are available.
The method assumes an LMO-compatible feasible region; in the context of quadratic programming, this means in particular that quadratic constraints have been transformed to objective penalties.
The heuristic consists in starting from any point in the feasible region, computed for instance from the LMO with an arbitrary direction.
From that point, the heuristic will compute the gradient, and call the LMO to compute another extreme point.
The algorithm iterates for a certain budget of maximum iterations or until it cycles to a vertex already encountered.
Follow-the-gradient can be viewed as a (nonconvergent) Frank-Wolfe algorithm with a unit step size.
Interestingly, the approximation guarantees of our follow-the-gradient heuristics matches the lower bounds of~\cite{baes2012minimizing};
see also~\cite{HPW2022} for a proximity result that provides good intuition why such a heuristic can be often powerful.
Preliminary computational experiments highlighted that gradient-based heuristics are not essential to the performance of our overall framework, therefore we do not report results for this heuristic in the main text.
The reason is likely that the solutions it obtains are redundant with the ones found during the execution of Frank-Wolfe.
This heuristic could, however, prove useful when using Boscia as an exact solver instead, with the advantage of only assuming access to the constraints as a linear minimization oracle, or also as seeding for large neighborhood search heuristics (see next section).

\subsection{Large neighborhood search}
Large Neighborhood Search (LNS) heuristics are powerful techniques for
mixed-integer optimization. Many well-known LNS heuristics for MIP rely on LP
relaxation solutions~\cite{RENS,RINS}, or reference solutions generated in other components such as cutting
planes~\cite{bolusani2024multi} or other primal heuristics~\cite{mexi2023using}.
In all these approaches, the reference solutions are fractional and feasible for the LP relaxation.
We propose instead to leverage the
current Frank--Wolfe iterate and the integer-feasible vertices sampled during the
execution of the Frank-Wolfe algorithm within the branch-and-bound process. Note that
these vertices are not necessarily feasible for the original MIQCQP continuous relaxation since 
we relax quadratic constraints via penalties.
\paragraph{Active Set Enforced Neighborhood Search.}
 Algorithm~\ref{alg:asens-generic} introduces
\textit{Active Set Enforced Neighborhood Search (ASENS)}, a new heuristic that
capitalizes on the structure of the FW active set by searching in the
neighborhood spanned by the active set points.
Recent FW variants maintain small active sets through
correction and local steps (see~\cite{tsuji2022pairwise} for the blended
pairwise conditional gradients variant
and~\cite{besancon2025improvedalgorithmsnovelapplications} for guarantees when
favoring local steps). We apply ASENS only when more than 50\% of the variables take
the same value across all active set vertices; ASENS then fixes these variables
accordingly and solves the resulting subproblem. In practice, the subproblem
often fixes substantially more than half of the variables, and we also restrict
the continuous variables to the convex hull of the active set. Preliminary
experiments showed that follow-the-gradient is often redundant with iterates
encountered during FW itself, while ASENS provides complementary improvements.

\begin{algorithm}[ht]
  \caption{Active Set Enforced Neighborhood Search (ASENS)}
  \label{alg:asens-generic}
  \begin{algorithmic}[1]
  \Require active set vertices $S=\{v^{1},\dots,v^{t}\}$, MIQCQP as in~\eqref{eq:miqcqp}, threshold $\tau \in (0,1)$
  \State Let $J=\{j:\;v^{1}_j = \dots = v^{t}_j\}$
  \If{$|J|/n < \tau$}
    \Comment{not enough variables fixed}
    \State \Return
  \EndIf
  \State Fix $x_j \gets v^{(1)}_j$ for all $j\in J$
  \State Restrict bounds of $x_j$ to $[\min_{i} v^{i}_j, \max_{i} v^{i}_j]$ for all $j\notin J$
  \State Solve the resulting MIQCQP
\end{algorithmic}
\end{algorithm}

\paragraph{Undercover.}
We implement a variant of the undercover heuristic~\cite{Undercover},
specifically designed for nonlinear problems. Undercover (UC) identifies and
fixes a subset of variables such that the resulting subproblem becomes linear
and, in many cases, easier to solve. Our implementation solves a minimum vertex
cover integer linear program to determine this subset, then fixes these
variables to values from the current Frank--Wolfe iterate from Boscia's
branch-and-bound tree. We also experimented with undercover initialized from
multiple diverse covers but did not observe performance improvements.

\paragraph{RINS.}
We implement a Frank--Wolfe-based variant of the Relaxation Induced Neighborhood
Search (RINS)~\cite{RINS}. While traditional RINS fixes variables that have
identical values in both the incumbent solution and the current relaxation
solution, our variant identifies promising search regions by fixing all
variables that have identical values in both the incumbent solution and the
current Frank--Wolfe iterate from the branch-and-bound tree.
Similar to ASENS, RINS is used only when more than 50\% of variables agree between the two
solutions.

\subsection{Specialized heuristics for QUBOs}
For QUBOs whose objective matrix has a \emph{bipartite} underlying graph, an initial solution can be improved on by alternately optimizing over each component of the graph.
This heuristic can seem redundant in comparison with the more general approach described above.
However, as there is a simple closed-form solution for each step of the alternating optimization, this specific case can be implemented very efficiently.
The inspiration for this comes from similar problems encountered in quantum communications, namely, the computation of the local bound of a (bipartite) Bell inequality~\cite{designolle2023improved}.
In practice, however, this approach was consistently outperformed by the other heuristics described above.

\section{Computational experiments}\label{sec:experiments}

In this section, we evaluate the performance of our heuristic and the impact of its different components.
We seek to answer the following questions:
\begin{itemize}
\item How does our heuristic perform on a diverse set of MIQCQPs?
\item What is the impact of parallelization on solution quality and efficiency?
\item How do the different LNS heuristics contribute to the overall performance?
\item What are the effects of varying the number of Frank-Wolfe iterations, power penalty parameter \(p\),
and convexification parameter \(\ell\)?
\end{itemize}

We measure performance with respect to different metrics, including the time to the first feasible solution (TTF),
the gap at the end of the time limit (compared to a reference dual bound), and the primal integral (PI).
These metrics highlight different important criteria that can be relevant depending on the context.

\subsection{Setup}

Our approach is implemented in Julia 1.11.4 on top of the Boscia and Frank-Wolfe packages.
As underlying MIP solver, we use Gurobi 12.0.1.
For our computational experiments, we used an Intel(R) Xeon(R) Gold 5122 @3.6GHz.
We used 8 threads and a memory limit of 96~GB RAM.
The time limit was set to 300 seconds. We test our method on the 319 instances of the QPLIB benchmark set~\cite{furini2019qplib} with discrete variables.
The instances are available at \url{https://qplib.zib.de}.

\paragraph{Evaluation metrics.}
To measure the performance of our heuristic, we employ standard metrics, such as the
time to first feasible solution, primal gap, and primal integral~\cite{berthold2015heuristic}.
The time to first feasible solution (TTF) is the time taken to find the first feasible solution during the search process.
The primal gap (Gap) of a heuristic solution $\tilde{x}$ with objective value $\tilde{f} = f(\tilde{x})$ with respect to the best known solution $x^*$ with objective value $f^* = f(x^*)$ is defined as
  \[
  \gamma(\tilde{x}) =
  \begin{cases}
    0 & \text{if } |\tilde{f}| = |f^*| = 0, \\
    1 & \text{if } \tilde{f} \cdot f^* < 0,\,\text{or no solution is found}, \\
    \dfrac{|\tilde{f} - f^*|}{\max(|\tilde{f}|, |f^*|)} & \text{otherwise}.
  \end{cases}
  \]
  The primal integral captures the evolution of the primal gap over time.
  Let $t_0=0$ be the start of the search process, $t_1, \ldots, t_{s-1}$ be the points in time when a new incumbent solution is found,
  and $t_s = T$ be the end of the search process, then the primal integral is defined as
  \[
  P(T) = \sum_{i=1}^{s} \gamma(\tilde{x}^i) \cdot (t_i - t_{i-1}).
  \]
  If the heuristic does not find any feasible solution, then the primal integral is equal to the time limit $T$.
  To report average values for these metrics across a set of instances, we use the shifted geometric mean with a shift of 1.

\paragraph{Parallelization strategy.}
Our parallel implementation strategy executes our algorithm using varying configurations in each thread. Specifically, we
ran our heuristics with different power penalty parameters $p$ ranging from
1.2 to 1.8. For problems with binary quadratic objectives, we employed various
convexification strategies to explore diverse solution approaches
simultaneously. The choice of different values for $p$ and convexification
parameters is motivated by preliminary experiments that indicate there is no
single best value for these parameters across all instances (see Section~\ref{sec:convexification_results}), and trying different
values allows us to explore different regions of the solution space.

To further diversify the search and avoid getting trapped in local optima, we
implemented a restart strategy that reinitializes our algorithm every 10 to 1000
nodes, depending on the problem characteristics. Each restart either uses the
current best solution as a warm start or initiates the search with a random
gradient direction to explore different regions of the solution space.

\section{Results}
\label{sec:results}

In Table~\ref{tab:base_7} we report our finding on QPLIB and the two categories of instances:
MIQCQP and MIQP (without quadratic constraints).
For each category we report the following metrics: the number of instances for which we found a
feasible solution, and for instances where we found a solution, we report
the average optimality gap compared to
the best known solution value from QPLIB, the primal integral also
using the best known solution as a reference, and time to the first feasible solution.

\begin{table}[ht]
\centering
\caption{Performance of the parallel heuristic on QPLIB instances.}
\label{tab:base_7}
  \setlength{\tabcolsep}{1.2em}
  \scriptsize
\begin{tabular}{lrrrr}
\toprule
Category & Found & Gap (\%) & PI & TTF\\
\midrule
MIQP & 163/166 & 4.06 & 5.23 & 1.76\\
MIQCQP & 98/153 & 25.30 & 50.56 & 7.31\\
\midrule
Total & 261/319 & 11.57 & 12.77 & 3.18\\
\bottomrule
\end{tabular}
\end{table}

Our heuristic is particularly effective for MIQP instances, and is able to
find feasible solutions for all but three instances in this category, and achieves an average
optimality gap of 4.06\% on instances where we found a solution.
Notably, we achieved optimality (0\% gap) on 104 instances.
The MIQCQP category is more challenging, and we are able to find
feasible solutions for 98 out of 153 instances, with an average optimality gap of 25.30\%.
We achieved optimality on 28 instances.
In total we found feasible solutions for 261 out of 319 instances, with an average
time to first feasible solution of 3.18 seconds. The average optimality gap was
11.57\%, with the average primal integral at 12.77.

\subsection{Competition and improvements on QPLIB}

Our heuristic achieved first place in the Land-Doig MIP Computational Competition 2025 \url{https://www.mixedinteger.org/2025/competition},
demonstrating its effectiveness on challenging MIQCQPs. Furthermore, the heuristic found eight new best-known solutions
for QPLIB instances, which are listed in Table~\ref{tab:qplib} and can be found at \url{https://qplib.zib.de/}.
We highlight that this improvement was achieved within the five-minute time limit of the competition.

\begin{table}[h]
  \centering
  \caption{Improved solutions for QPLIB instances. The columns indicate the
  instance name, the objective sense, the previous best-known solution value,
  the new solution value found by our heuristic, and the relative improvement
  (gap) in percentage.}
  \label{tab:qplib}
  \scriptsize
  \begin{tabular}{lrrrr}
    \toprule
    Instance & Obj.~Sense & Previous Best & New Solution & Gap (\%)\\
    \midrule
    QPLIB\_2169 & max & 29.0 & 30.0 & 3.45 \\
    QPLIB\_2174 & max & 150.0 & 152.0 & 1.33 \\
    QPLIB\_2205 & max & 88.0 & 90.0 & 2.27 \\
    QPLIB\_3347 & min & 3,819,920 & 3,818,879 & 0.03 \\
    QPLIB\_3584 & min & -25,254 & -25,386 & 0.52 \\
    QPLIB\_3709 & min & 5,726,530 & 5,710,645 & 0.28 \\
    QPLIB\_3860 & min & -19,685 & -20,161 & 2.42 \\
    QPLIB\_10022 & min & 6,267,782 & 1,374,066 & 78.07 \\
    \bottomrule
  \end{tabular}
\end{table}

\subsection{Comparison to state-of-the-art solvers} 

To assess the quality of our heuristic relative to state-of-the-art solvers, we compare its performance against Gurobi 12.0.1 and the open-source solver SCIP (pre-release version 10.0.0, GitHub commit \texttt{30844abadc}\footnote{\url{https://github.com/scipopt/scip}})~\cite{bolusani2024scip} on the QPLIB benchmark instances.
Both solvers were run with default settings and a time limit of 300 seconds on the same hardware as our heuristic.
Table~\ref{tab:comparison} summarizes the results, comparing the solution quality achieved by each solver within the time limit.

\begin{table}[ht]
  \centering
  \caption{Performance comparison of the parallel heuristic against Gurobi and SCIP on QPLIB instances with a 300-second time limit.
  The table shows the number of instances where the heuristic found a solution with the same, better, or worse objective value compared to each solver.
  The ``None'' column indicates instances where no solution was found by either method.}
  \label{tab:comparison}
  \renewcommand{\arraystretch}{1.15}
  \setlength{\tabcolsep}{0.6em}
  \small
  \begin{tabular}{@{}l*{8}{c}@{}}
    \toprule
    & \multicolumn{4}{c}{vs.\ Gurobi} & \multicolumn{4}{c}{vs.\ SCIP} \\
    \cmidrule(lr){2-5} \cmidrule(lr){6-9}
    Type & Same & Better & Worse & None & Same & Better & Worse & None \\
    \midrule
    MIQP   & 93 & 30 & 43 &  0 & 40 & 94 & 31 &  1 \\
    MIQCQP & 27 &  1 & 101 & 24 & 15 & 26 & 81 & 31 \\
    \bottomrule
  \end{tabular}
\end{table}

Our heuristic demonstrates competitive performance, particularly on MIQP instances, where it finds better solutions than Gurobi on 30 instances and better solutions than SCIP on 94 instances.
For MIQCQPs, both Gurobi and SCIP outperform our heuristic, highlighting the substantial challenges posed by quadratic constraints.

\subsection{Impact of parallelization}
Next, we compare the performance of our heuristic when run with a single master thread
(base) versus when run with multiple threads (heur-parallel), with each thread exploring a different part of the search space.
For a fair comparison, even when a single master thread is used, our heuristic continues to
utilize all available threads for the LMO.
Table~\ref{tab:base_base_7} summarizes the results of this comparison over instances solved by at least one of the two settings,
and over instances solved by both settings.

As shown in Table~\ref{tab:base_base_7}, distributing the search across multiple threads yields
greater benefits than concentrating resources on the LMO calls.
Specifically, the multithreaded version of our approach solves 8 more instances than the
single-threaded version, reducing the average optimality gap from 16.89\% to 12.07\% and
decreasing the average primal integral from 17.40 to 13.10, considering all instances solved
by at least one of the two settings.

\begin{figure}[ht]
  \centering
  \caption{
    Distribution of the optimality gap, and primal integral over MIQCQPs where a solution was found by at least one of the settings.
  }
  \includegraphics[width=\textwidth]{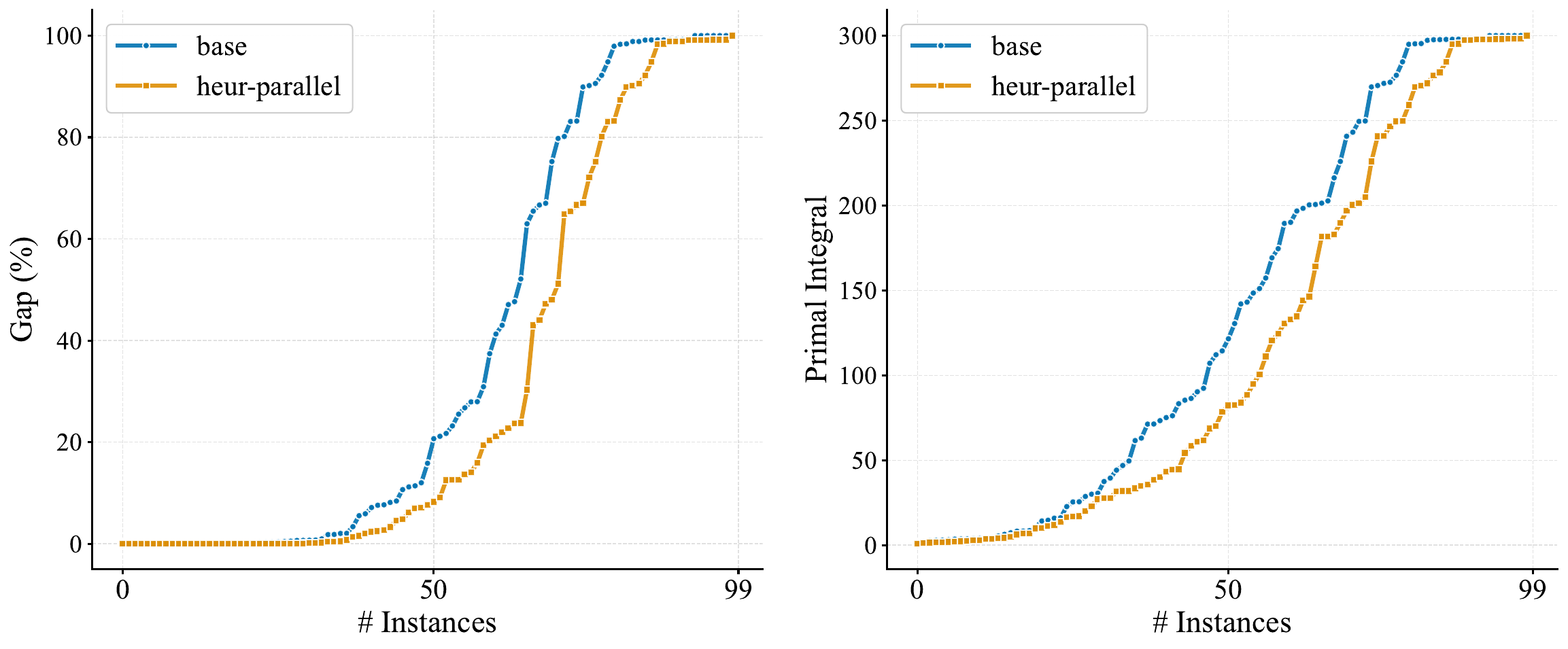}
  \label{fig:time_first_base_7_miqcqp}
\end{figure}

Over MIQPs, there is no significant difference in the number of instances solved, or the average performance metrics.
However, for MIQCQPs, the multithreaded version solves 6 more instances and significantly reduces the average optimality gap and primal integral.
Figure~\ref{fig:time_first_base_7_miqcqp} shows the distribution of the optimality gap, and primal integral over MIQCQPs for both settings.

\begin{table}[ht]
\centering
\caption{
  Performance comparison between single-threaded (base) and multi-threaded (heur-parallel) approaches on QPLIB instances.
  Bold values are the best ones: highest number of instances solved, lowest averaged optimality gap, primal integral, and time to first feasible solution.
}
\label{tab:base_base_7}
\setlength{\tabcolsep}{0.6em}
\scriptsize
\begin{tabular}{@{}llrrrrrrr@{}}
\toprule
Type & Setting & \multicolumn{3}{c}{At Least One Solved} & \multicolumn{4}{c}{All Solved} \\
\cmidrule(lr){3-5} \cmidrule(lr){6-9}
 & & Found & Gap (\%) & PI & Found & Gap (\%) & PI & TTF \\
\midrule
\multirow{2}{*}{All (263)} & base & 256 & 14.10 & 16.02 & 254 & 12.46 & 14.55 & 2.94 \\
 & heur-parallel & \win{261} & \win{12.07} & \win{13.10} & 254 & \win{11.56} & \win{12.02} & \win{2.90} \\
\midrule
\multirow{2}{*}{MIQP (164)} & base & \win{164} & 4.51 & 6.40 & 163 & 4.53 & 6.32 & \win{1.74} \\
 & heur-parallel & 163 & \win{4.47} & \win{5.38} & 163 & \win{4.06} & \win{5.23} & 1.76 \\
\midrule
\multirow{2}{*}{MIQCQP (99)} & base & 92 & 31.97 & 66.75 & 91 & 28.20 & 58.91 & 6.54 \\
 & heur-parallel & \win{98} & \win{25.89} & \win{51.49} & 91 & \win{26.37} & \win{47.73} & \win{6.24} \\
\bottomrule
\end{tabular}
\end{table}

\subsection{Impact of LNS heuristics}

To assess the impact of the LNS heuristics finding feasible solutions, we compare the performance of the
base heuristic with the LNS-off variant, which does not use any LNS heuristics,
the ASENS-only variant, which only uses the ASENS heuristic, and the UC-only variant, which only uses the undercover heuristic.
All variants still use the rounding heuristics described in Section~\ref{sec:rounding_heuristics}.
The results are summarized in Table~\ref{tab:base_heur_prob_00_heur_only_arens_heur_only_undercover} and demonstrate both the
significant contribution of the LNS heuristics to our framework's performance and the robustness of the underlying approach without them.

\begin{figure}[ht]
  \centering
  \caption{
    Distribution of the optimality gap, and primal integral over MIQCQPs where a solution was found by at least one of the settings.
  }
  \includegraphics[width=\textwidth]{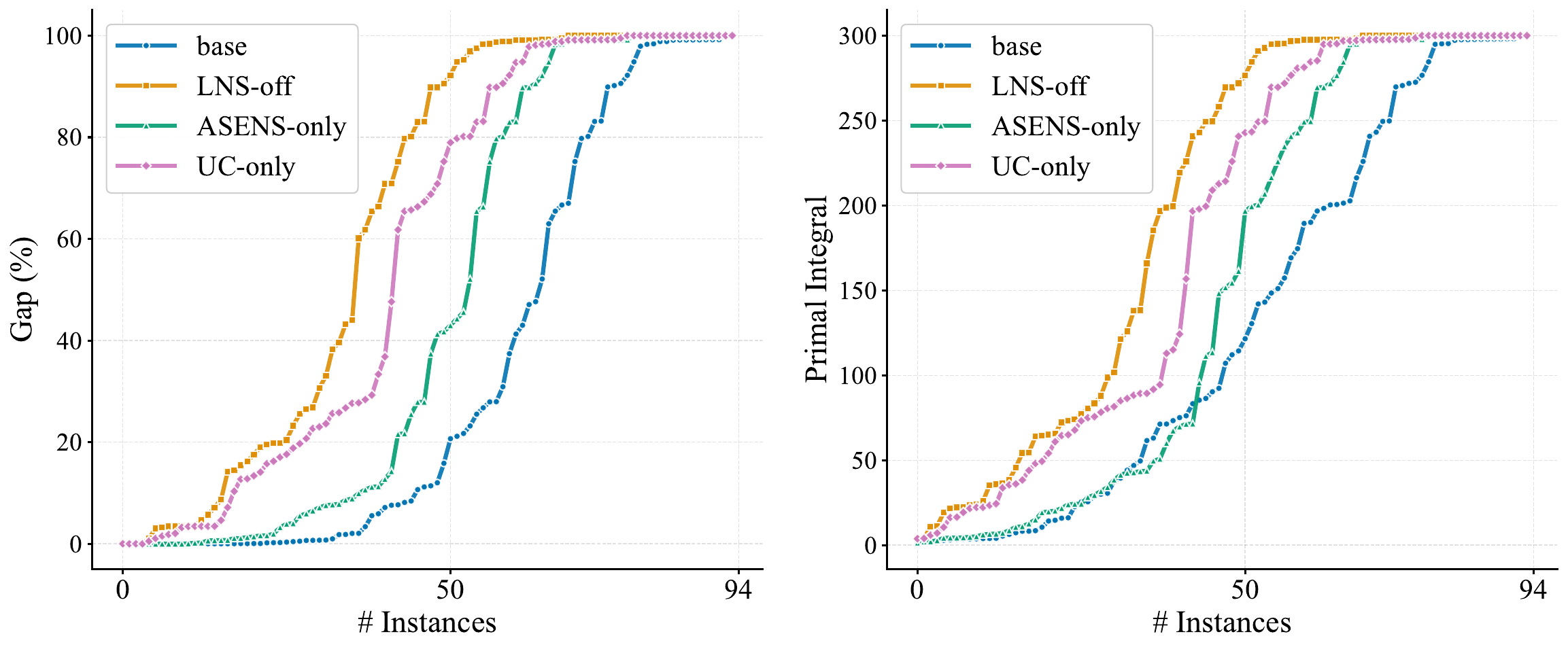}
  \label{fig:time_first_base_heur_prob_00_heur_only_arens_heur_only_undercover}
\end{figure}

In particular, the LNS-off setting, which relies exclusively on LMO calls without any
quadratic subproblem solving, successfully finds feasible solutions for a
substantial amount of 226 instances, with an average optimality gap of 28.67\% and a primal integral of 29.42
over the 258 instances where at least one of the settings found a solution.
This result demonstrates that our framework's core mechanism of leveraging Frank-Wolfe
iterations over the convex hull of integer-feasible points is robust even when
restricted to LMO calls.
Similarly, the UC-only setting performs well, finding feasible solutions for 249 instances, also without
solving any quadratic subproblems.
The ASENS-only setting leads to an increased number of feasible solutions found, when compared to LNS-off,
and achieves the closest performance with respect to the primal integral and optimality gap to the base heuristic.
As shown in Table~\ref{tab:base_heur_prob_00_heur_only_arens_heur_only_undercover} and in
Figure~\ref{fig:time_first_base_heur_prob_00_heur_only_arens_heur_only_undercover}, the contribution of the LNS heuristics is particularly evident in the MIQCQP category, where the base heuristic finds feasible solutions for 92 instances,
while the LNS-off setting finds solutions for only 74 instances, and the ASENS-only and UC-only settings find solutions for 86 and 85 instances, respectively.

\begin{table}[ht]
\centering
\caption{
  Impact of different LNS heuristics on the performance across QPLIB instances.
  Bold values are the best ones: highest number of instances solved, lowest averaged optimality gap, primal integral, and time to first feasible solution.
}
\label{tab:base_heur_prob_00_heur_only_arens_heur_only_undercover}
\setlength{\tabcolsep}{0.6em}
\scriptsize
\begin{tabular}{@{}llrrrrrrr@{}}
\toprule
Type & Setting & \multicolumn{3}{c}{At Least One Solved} & \multicolumn{4}{c}{All Solved} \\
\cmidrule(lr){3-5} \cmidrule(lr){6-9}
 & & Found & Gap (\%) & PI & Found & Gap (\%) & PI & TTF \\
\midrule
\multirow{4}{*}{All (258)} & base & \win{256} & \win{12.86} & \win{15.10} & 225 & \win{12.29} & \win{12.29} & 2.24 \\
 & LNS-off & 228 & 28.67 & 29.42 & 225 & 21.22 & 21.26 & 2.13 \\
 & ASENS-only & 243 & 18.66 & 16.53 & 225 & 14.12 & 12.21 & \win{2.02} \\
 & UC-only & 249 & 24.34 & 26.44 & 225 & 20.44 & 20.50 & 2.05 \\
\midrule
\multirow{4}{*}{MIQP (164)} & base & \win{164} & \win{4.51} & \win{6.40} & 153 & \win{3.65} & \win{5.29} & 1.44 \\
 & LNS-off & 154 & 14.56 & 11.82 & 153 & 10.56 & 9.42 & 1.45 \\
 & ASENS-only & 157 & 7.79 & 6.64 & 153 & 4.56 & 5.14 & 1.43 \\
 & UC-only & \win{164} & 10.97 & 10.74 & 153 & 10.46 & 9.23 & \win{1.39} \\
\midrule
\multirow{4}{*}{MIQCQP (94)} & base & \win{92} & \win{29.08} & \win{61.58} & 72 & \win{33.10} & \win{64.23} & 4.93\\
 & LNS-off & 74 & 57.92 & 138.65 & 72 & 47.40 & 110.62 & 4.27 \\
 & ASENS-only & 86 & 40.16 & 72.93 & 72 & 37.46 & 66.32 & \win{3.83} \\
 & UC-only & 85 & 51.65 & 119.61 & 72 & 44.72 & 103.18 & 4.12 \\
\bottomrule
\end{tabular}
\end{table}

\begin{figure}[ht]
  \centering
  \caption{
    Impact of the RINS heuristic on the performance of the base heuristic across MIQCQPs where a solution was found by at least one of the settings.
  }
  \includegraphics[width=\textwidth]{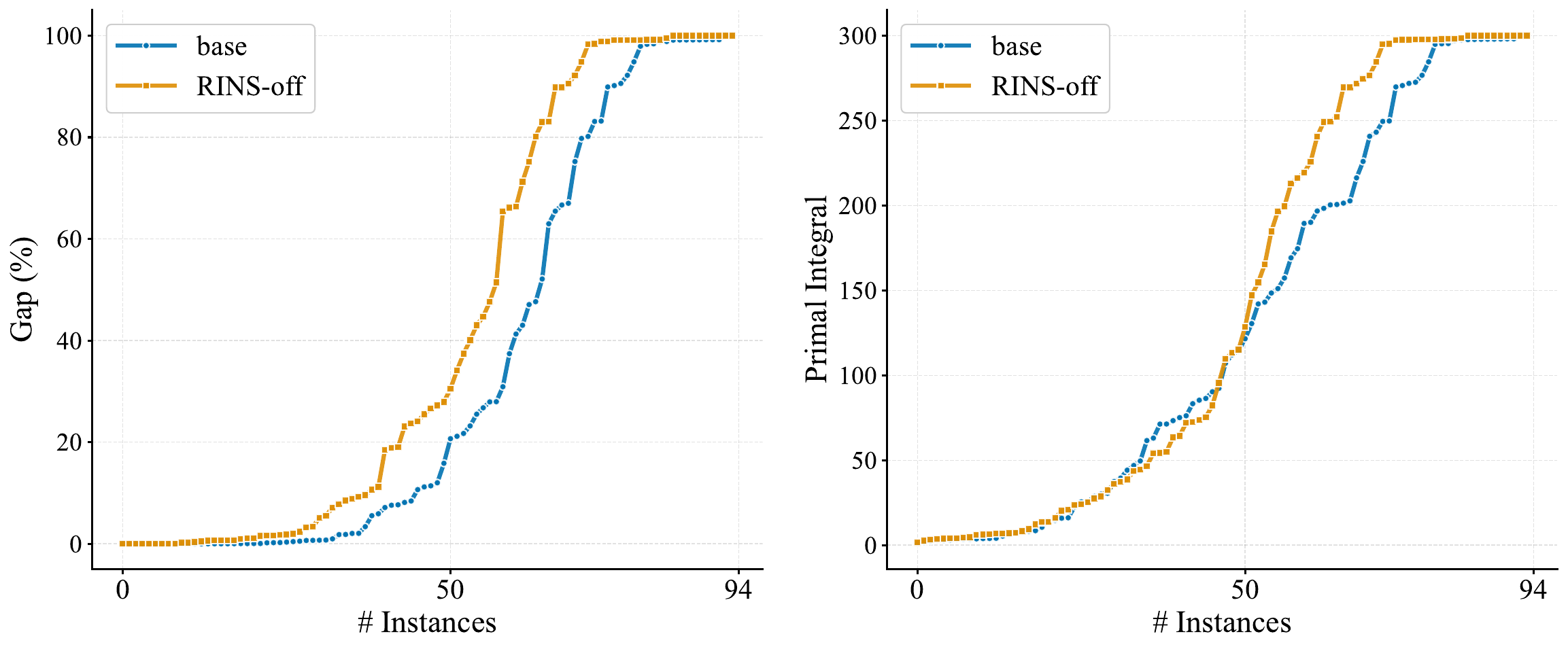}
  \label{fig:rins_impact_miqcqp}
\end{figure}

Finally, our base heuristic also includes RINS. To measure its impact, we compare
the base heuristic (with RINS) against a RINS-off variant. Note
that RINS, as an improvement heuristic, focuses on refining existing solutions.
The impact of RINS on MIQPs is moderate since the base heuristic already finds
good quality solutions for most instances. In particular, removing RINS increases the
average gap and primal integral by only around 0.5\%. However, for MIQCQPs,
RINS provides improvements in solution quality as shown in
Figure~\ref{fig:rins_impact_miqcqp}. Specifically, including RINS reduces the
average optimality gap by 5.8\% and decreases the average primal integral by a value of 11.5.

\subsection{Impact of the Frank-Wolfe accuracy}

We now assess the impact of the Frank-Wolfe accuracy on the performance of our heuristic, specifically by varying the number of maximum Frank-Wolfe iterations.
The results for this experiment are summarized in Table~\ref{tab:fw_iter_001_fw_iter_010_base_fw_iter_1000}.

\begin{table}[ht]
\centering
\caption{
  Performance of the base heuristic with different limits on the number of Frank-Wolfe iterations.
  Bold values are the best ones: highest number of instances solved, lowest averaged optimality gap, primal integral, and time to first feasible solution.
}
\label{tab:fw_iter_001_fw_iter_010_base_fw_iter_1000}
\setlength{\tabcolsep}{0.6em}
\scriptsize
\begin{tabular}{@{}llrrrrrrr@{}}
\toprule
Type & Setting & \multicolumn{3}{c}{At Least One Solved} & \multicolumn{4}{c}{All Solved} \\
\cmidrule(lr){3-5} \cmidrule(lr){6-9}
 & & Found & Gap (\%) & PI & Found & Gap (\%) & PI & TTF \\
\midrule
\multirow{4}{*}{All (262)} & fw-iter=1 & \win{257} & 17.80 & 24.63 & 249 & 16.09 & 22.23 & 2.79 \\
 & fw-iter=10 & 256 & \win{13.85} & \win{15.84} & 249 & \win{12.64} & \win{14.12} & \win{2.74}  \\
 & fw-iter=100 & 253 & 16.65 & 17.20 & 249 & 14.49 & 15.26 & 2.89 \\
 & fw-iter=1000 & 253 & 17.50 & 19.78 & 249 & 15.40 & 17.52 & 3.13 \\
\midrule
\multirow{4}{*}{MIQP (164)} & fw-iter=1 & 164 & 8.22 & 12.79 & 164 & 8.22 & 12.79 & 1.81 \\
 & fw-iter=10 & 164 & \win{4.51} & \win{6.40} & 164 & \win{4.51} & \win{6.40} & \win{1.79} \\
 & fw-iter=100 & 164 & 6.90 & 7.21 & 164 & 6.90 & 7.21 & 1.80 \\
 & fw-iter=1000 & 164 & 7.99 & 8.12 & 164 & 7.99 & 8.12 & 1.85 \\
\midrule
\multirow{4}{*}{MIQCQP (98)} & fw-iter=1 & \win{93} & 35.77 & 71.33 & 85 & 32.93 & 62.54 & 5.75 \\
 & fw-iter=10 & 92 & \win{31.41} & \win{65.72} & 85 & \win{30.17} & \win{59.04} & \win{5.58} \\
 & fw-iter=100 & 89 & 35.00 & 68.00 & 85 & 30.70 & 59.82 & 6.32 \\
 & fw-iter=1000 & 89 & 35.34 & 81.45 & 85 & 31.16 & 71.63 & 7.46 \\
\bottomrule
\end{tabular}
\end{table}

We aim to strike a balance between computational effort in the convex relaxations and solution accuracy.
Too few iterations would cause the relaxation solutions not to move sufficiently from one node to the next, while too many iterations would make the framework spend too much time computing high-accuracy solutions, which will in any case only serve as starting point to compute integer-feasible solutions through, e.g., rounding afterwards.
Our experiments show that 10 to 100 FW iterations per node yields a good balance between accuracy and relaxation cost, in terms of number of instances with at least one solution, mean residual gap, and primal integral.

\subsection{Impact of convexification and power penalty reformulation parameters.}
\label{sec:convexification_results}
For the 124 MIQPs with binary variables, we tested various convexification parameters.
Table~\ref{tab:base-convexify-objective-0-6-convexify-objective-0-7-convexify-objective-0-9-convexify-objective-1-0}
shows that there is no single best convexification parameter for all instances.
Interestingly, fully convexifying the objective function leads to a significant increase in the average gap and primal integral.
\begin{table}[ht]
\centering
\caption{
  Performance of the base heuristic with different convexification parameters $\ell$.
  Bold values are the best ones: highest number of instances solved, lowest averaged optimality gap, primal integral, and time to first feasible solution.
}
\label{tab:base-convexify-objective-0-6-convexify-objective-0-7-convexify-objective-0-9-convexify-objective-1-0}
\setlength{\tabcolsep}{0.6em}
\scriptsize
\begin{tabular}{@{}llrrrr@{}}
\toprule
Type & Setting & \multicolumn{4}{c}{All Solved} \\
\cmidrule(lr){3-6}
 & & Found & Gap (\%) & PI & TTF \\
\midrule
\multirow{5}{*}{Binary QP (124)} & $\ell=0.6$ & 124 & 3.66 & 5.17 & 1.46 \\
 & $\ell=0.7$ & 124 & 3.35 & 5.09 & 1.45 \\
 & $\ell=0.8$ & 124 & \win{3.18} & \win{4.35} & 1.42 \\
 & $\ell=0.9$ & 124 & \win{3.18} & 4.79 & 1.46 \\
 & $\ell=1.0$ & 124 & 5.02 & 6.62 & \win{1.41} \\
\bottomrule
\end{tabular}
\end{table}

Further we evaluated the impact of the power penalty parameter \(p\) on the performance of the heuristic.
The results are summarized in Table~\ref{tab:base-power-penalty-1-2-power-penalty-1-3-power-penalty-1-4-power-penalty-1-6-power-penalty-1-7-power-penalty-1-8}.
We observe that the performance of the heuristic is not significantly affected by the choice of \(p\) either.

\begin{table}[ht]
\centering
\caption{
  Performance of the base heuristic with different power penalty parameters $p$.
  Bold values are the best ones: highest number of instances solved, lowest averaged optimality gap, primal integral, and time to first feasible solution.
}
\label{tab:base-power-penalty-1-2-power-penalty-1-3-power-penalty-1-4-power-penalty-1-6-power-penalty-1-7-power-penalty-1-8}
\setlength{\tabcolsep}{0.6em}
\scriptsize
\begin{tabular}{@{}llrrrrrrr@{}}
\toprule
Type & Setting & \multicolumn{3}{c}{At Least One Solved} & \multicolumn{4}{c}{All Solved} \\
\cmidrule(lr){3-5} \cmidrule(lr){6-9}
 & & Found & Gap (\%) & PI & Found & Gap (\%) & PI & TTF \\
\midrule
\multirow{7}{*}{MIQCQP (96)} & p=1.2 & 88 & 33.59 & 65.19 & 82 & \win{30.67} & \win{53.42} & \win{4.93} \\
 & p=1.3 & 90 & 32.58 & 66.52 & 82 & 30.91 & 55.78 & 5.23 \\
 & p=1.4 & 86 & 35.32 & 68.04 & 82 & 29.82 & 55.29 & 5.12 \\
 & p=1.5 & 89 & 33.90 & 65.92 & 82 & 30.81 & 56.38 & 5.52 \\
 & p=1.6 & \win{91} & \win{31.34} & \win{64.63} & 82 & 30.82 & 57.50 & 5.14 \\
 & p=1.7 & 88 & 33.72 & 68.70 & 82 & 31.32 & 57.07 & 5.11 \\
 & p=1.8 & 86 & 37.01 & 69.11 & 82 & 31.18 & 55.60 & 5.00 \\
\bottomrule
\end{tabular}
\end{table}

These observations motivate our choice of parallel and restart strategies: diversifying the search to different values of $p$ and $\ell$ maximizes the chance of finding values suited to each instance.
We note that there should exist criteria to enlighten this choice based on the instance characteristics, but we leave this for further research.

\section{Conclusion}

In this paper, we presented a primal heuristic framework for solving
mixed-integer quadratically constrained quadratic programs. Our
approach builds upon the Frank-Wolfe-based branch-and-bound framework Boscia and
is designed to efficiently explore the solution space through gradient-guided
directions and large neighborhood search heuristics that exploit
integer-feasible vertices sampled during Frank-Wolfe iterations.
The framework achieved first place in the Land-Doig MIP Computational Competition 2025 and
discovered eight new best-known solutions for QPLIB instances within the
five-minute time limit of the competition, demonstrating its practical effectiveness
on challenging nonconvex mixed-integer optimization problems.

\section*{Acknowledgements}
We thank the organizing committee of the Land-Doig MIP Computational
Competition, with special thanks to Aleksandr Kazachkov and Jan Kronqvist.
Research reported in this paper was partially supported through
the Research Campus Modal funded by the German Federal Ministry of Research, Technology, and Space
(fund numbers 05M14ZAM, 05M20ZBM) and the Deutsche
Forschungsgemeinschaft (DFG) through the DFG Cluster of Excellence MATH+.
The work of Mathieu Besançon benefited from the support of the FMJH Program PGMO and
from MIAI at Université Grenoble Alpes (grant ANR-19-P3IA-0003).

\bibliographystyle{splncs04}
\bibliography{references}

\end{document}